\documentclass[11pt,amssymb]{article}
\usepackage{amsfonts,amsmath,latexsym}

\title{An end-to-end construction for compact constant mean curvature
surfaces}

\author{M. Jleli and F. Pacard
\thanks{LAMA, Universit\'e Paris 12, 61 Avenue du G\'en\'eral de Gaulle,
94010 Cr\'eteil, France} \\ Universit\'e Paris 12}

\date{}

\newtheorem{theorem}{Theorem}
\newtheorem{proposition}{Proposition}

\newtheorem{lemma}{Lemma}
\newtheorem{definition}{Definition}

\newcommand{\N}{\mathbb N}
\newcommand{\R}{\mathbb R}
\newcommand{\del}{\partial}
\newcommand{\calC}{{\mathcal C}}
\newcommand{\calL}{{\mathcal L}}
\newcommand{\calM}{{\mathcal M}}
\newcommand{\calE}{{\mathcal E}}
\newcommand{\calD}{{\mathcal D}}
\newcommand{\calW}{{\mathcal W}}
\newcommand{\calN}{{\mathcal N}}
\newcommand{\calK}{{\mathcal K}}

\begin{document}

\maketitle

\section{Introduction}

The theory of constant mean curvature surfaces in Euclidean space
has been the object of intensive study in the past years. In the
case of complete noncompact constant mean curvature surfaces, the
moduli space of such surfaces is now fairly well understood (at
least in the genus $0$ case) \cite{Kus-Maz-Pol},
\cite{Kar-Kus-Sul-1}, \cite{Kar-Kus-Sul-2} and many technics have
been developed to produce examples of such surfaces \cite{Kar},
\cite{Kap-1}, \cite{Maz-Pac}, \cite{Maz-Pac-Pol-2}.

\medskip

By contrast, the set of compact constant mean curvature surfaces
is not so well understood. In the early $80$'s, H. Wente has
constructed the first examples of genus $1$ constant mean
curvature surfaces \cite{Wen}. These genus $1$ surfaces have then
been thoughtfully studied by U. Pinkall and I. Sterling
\cite{Pin-Ste}. Examples of compact constant mean curvature
surface of higher genus are due to N. Kapouleas. In the genus $2$
case \cite{Kap-3}, these surfaces are obtained by "fusing" Wente
tori while in the case where the genus is greater than or equal to
$3$, these surfaces are obtain by connecting together large number
of mutually tangent unit spheres, using small catenoid necks
\cite{Kap-2}.

\medskip

In this paper, we would like to explain how the current knowledge
on the set of complete noncompact constant mean curvature surfaces
can be exploited to produce new examples of compact constant mean
curvature surfaces of genus greater than or equal to $3$.

\medskip

Our construction is based on three important tools which have been
developed for the understanding of complete noncompact constant
mean curvature surfaces~:
\begin{itemize}
\item[(i)] The moduli space theory as developed by R. Kusner, R.
Mazzeo and D. Pollack \cite{Kus-Maz-Pol}.

\item[(ii)] The end-addition result which has been developed by R.
Mazzeo, F. Pacard and D. Pollack \cite{Maz-Pac-Pol-1},
\cite{Maz-Pac-Pol-2} to produce complete noncompact constant mean
curvature surfaces with prescribed ends.

\item[(iii)] The end-to-end construction which was developed by J.
Ratzkin \cite{Rat} to connect two constant mean curvature surfaces
along their ends.
\end{itemize}

This ideas behind our construction can be described as follow~:
One can use the end-addition theory developed in
\cite{Maz-Pac-Pol-1}, \cite{Maz-Pac-Pol-2} to produce complete
constant mean curvature surfaces with prescribed Delaunay type
ends. This addition of ends procedure is quite flexible and one
can arrange so that the ends of these surfaces can be "plugged"
together to produce sequences (indexed by a discreet parameter $n
\in {\N}$) of compact surfaces which have mean curvature one
except in finitely many annular regions where their mean curvature
can be estimated by $1 + {\cal O} (e^{-\gamma n})$ for some
$\gamma > 0$. Next, one studies the mapping properties of the
Jacobi operator about this (almost) constant mean curvature
surface. To perform this analysis, we rely on the fact that
parametrices for the Jacobi operators on each complete noncompact
summand have been obtained in \cite{Kus-Maz-Pol} and we explain
how these can be glued together. Finally, it will remain to use a
standard perturbation argument to produce sequences of compact
constant mean curvature surfaces of arbitrary genus, greater than
or equal to $3$.

\medskip

The main advantage of our construction versus the one developed by
N. Kapouleas is that it is technically simple (once the above
mentioned technics are understood !), paralleling the fact that
the end-to-end construction of J. Ratzkin is simpler than the
previous constructions of complete noncompact surfaces. We obtain
a very precise description of the surfaces we produce. In
particular, our construction sheds light on the structure of the
set of compact constant mean curvature surfaces, showing that
these surfaces are isolated. Though this is probably a minor
point, the example of compact constant mean curvature surfaces we
obtain are geometrically different from the one obtained by N.
Kapouleas and cannot be obtained using his technic (roughly
speaking all the surfaces constructed by N. Kapouleas have close
to sequences of unit spheres linked by small catenoids and hence
have small injectivity radius while our examples do not
necessarily have small necks and hence have injectivity radius
uniformly bonded from below).

\medskip

Maybe a more important issue is the fact that our construction
points out interesting directions toward which the theory of
complete noncompact constant mean curvature surfaces should be
developed. Some properties of complete noncompact constant mean
curvature surfaces have been neglected and turn out to be
extremely important. This is for example the case of the notion of
"regular end" of a constant mean curvature surface (which is also
important in the construction of J. Ratzkin).

\medskip

Although our method can be applied to produce non symmetric
surface, the complete description of the set of compact surfaces
is far beyond our understanding, this is the reason why we have
chosen not to look for the utmost generality but to focuss on the
construction of symmetric surfaces. In order to explain the ideas
in our construction (keeping the technicalities as low as possible
and the notations as simple as possible), we will construct
constant mean curvature surfaces of arbitrary genus ($\geq 3$)
which have dihedral symmetry.

\medskip

Final remark, our construction generalizes to any dimension
\cite{Jle}.

\section{Delaunay surfaces}

In this section we recall some well known facts concerning the
family of Delaunay surfaces $D_\tau$ which are rotationally
invariant constant mean curvature surfaces in ${\mathbb R}^3$
\cite{Del}. We refer to \cite{Maz-Pac} for further details.

\subsection{Isothermal parametrization}

Delaunay surfaces can be parameterized, in isothermal coordinates,
by
\begin{equation}
X_\tau(s,\theta) =  \frac12\,\left(\tau\,e^{\sigma(s)}
\, \cos\theta,\tau\,e^{\sigma(s)}\,  \sin\theta,\kappa (s)\right),
\label{eq:2.1}
\end{equation}
for $(s, \theta) \in {\R}\times S^1$, where the functions $\sigma
$ and $\kappa $ are described as follows~: For any $\tau \in
(0,1]$, the function $\sigma $ is defined to be the unique smooth
nonconstant solution of the ordinary differential equation
\[
(\del_s\sigma)^2 + \tau^2 \cosh^2 \sigma =1, \qquad \del_s \sigma(0)
= 0,  \qquad \sigma(0) < 0,
\]
while, for any $\tau \in (-\infty,0)$, the function $\sigma $ is
defined to be the unique smooth nonconstant solution of the
ordinary differential equation
\[
(\del_s\sigma)^2 + \tau^2 \sinh^2 \sigma =1, \qquad \del_s \sigma(0)
= 0,  \qquad \sigma(0) <0 .
\]
Again, the definition of $\kappa $ differs according to the sign
of $\tau$. When $\tau \in (0,1]$, we define the function $\kappa$
by
\[
\del_s\kappa =\tau^2\, e^{\sigma}\, \cosh\sigma, \qquad \kappa(0) = 0,
\]
while when $\tau < 0$, we define the function $\kappa$ by
\[
\del_s \kappa =  \tau^2 \, e^{\sigma}\,\sinh \sigma, \qquad \kappa(0)
= 0.
\]
Observe that when $\tau>0$, $\kappa$ is monotone increasing,  and
hence $X_\tau$ is an embedding, whereas when $\tau < 0$, this is
no longer true and the surfaces are only immersed. The embedded
(resp. immersed) Delaunay surfaces $D_\tau$ are known as unduloids
(resp. nodoids). The parameter $\tau$ will be refered to as the
Delaunay parameter.

\medskip

As noted above, these surfaces are all periodic because the
functions $\sigma$ are. When $\tau \neq 1$, we define $s_\tau$ to
be equal to the {\em half} of the least period of $\sigma$. The
physical least period of the Delaunay surface $D_\tau$ is given by
\[
2 \, T_\tau : =  \frac{1}{2} \, \kappa_\tau (2 \, s_\tau)
\]

\noindent
{\bf Warning :} We agree that $2 \, s_\tau$ is equal to
the least period of $\sigma$ and $2\, T_\tau$ is the least period
of $D_\tau$.

\medskip

We claim that
\begin{lemma}
For all $\tau \in (-\infty , 0) \cup (0, 1]$, we have $\del_\tau
T_\tau > 0$.
\end{lemma}
{\bf Proof :} Observe that $\del_s \sigma >0$ on $(0, s_\tau)$ and
$\del_s \sigma <0$ on $(s_\tau , 2 \, s_\tau)$. Hence, for $s \in
(0, s_\tau)$, we can use $\sigma$ as a change of variable and
express $\kappa$ as a function of $\sigma \in (- \sigma_* ,
\sigma_*)$ where $\sigma_* >0$ satisfies $\tau^2 \, \cosh^2
\sigma_* =1$ when $\tau \in (0, 1]$ and $\sigma_* >0$ satisfies
$\tau^2 \, \sinh^2 \sigma_* =1$ when $\tau <0$.

\medskip

When $\tau <0$, we get
\[
2 \, T_\tau = \int_{- \sigma_*}^{\sigma_*} \frac{\tau^2 \,
e^\sigma \, \sinh \sigma}{\sqrt{1 -\tau^2 \sinh^2 \sigma} }\,
d\sigma =\int_{- \sigma_*}^{ \sigma_*} \frac{\tau^2 \, \sinh^2
\sigma }{\sqrt{1 -\tau^2 \sinh^2 \sigma}} \, d\sigma
\]
Performing the change of variable $\tau \, \sinh \sigma = \sin x$
we conclude that
\[
2 \, T_\tau = \int_{-\pi/2}^{\pi/2} \frac{\sin^2 x}{\sqrt{\tau^2
+ \sin^2 x}} \, dx
\]
which clearly implies that $\del_\tau T_\tau > 0$ when $\tau <0$.

\medskip

When $\tau >0$, we have
\[
2 \, T_\tau = \int_{-\sigma_*}^{\sigma_*} \frac{\tau^2 \, e^\sigma
\, \cosh \sigma}{\sqrt{1 -\tau^2 \cosh^2 \sigma} } \, d\sigma =
\int_{- \sigma_*}^{\sigma_*} \frac{\tau^2 \, \cosh^2 \sigma
}{\sqrt{1 -\tau^2 \cosh^2 \sigma} } \, d\sigma
\]
Performing the change of variable $\tau \, \sinh \sigma =
\sqrt{1-\tau^2} \, \sin x$ we conclude that
\[
2 \, T_\tau = \int_{- \pi/2}^{\pi/2} \sqrt{1 - (1-\tau^2) \,
\cos^2 x} \, dx
\]
which again implies that $\del_\tau T_\tau > 0$.
\hfill $\Box$

\medskip

It will be convenient to denote by $D_\tau^{\vec{a}}$ (resp.
$X_\tau^{\vec{a}}$) the Delaunay surface (resp. the
parameterization of the Delaunay surface) whose Delaunay parameter
is equal to $\tau$, whose axis is the line directed by $\vec{a}$
passing through the origin, and which has a neck passing through
the plane $\vec{x} \cdot \vec{a} = 0$. In particular, $D^{\vec
a}_\tau$ is invariant under the symmetry with respect to the plane
whose normal is $\vec a$. Granted this notation, the Delaunay
surface $D_\tau$ described (\ref{eq:2.1}) is equal to $D^{\vec
e_3}_\tau$ and $D^{\vec{a}}_\tau$ is obtained from $D_\tau^{\vec
e_3}$ by applying a rotation which sends $\vec e_3$ to $\vec a$.
Given a vector $\vec b$, the surface $D_\tau^{\vec a} + \vec b$ is
the surface $D_\tau^{\vec a}$ which has been translated by $\vec
b$. It is parameterized by $X_\tau^{\vec a} + \vec b$ and, when
$\vec b=0$, the unit normal vector field compatible with the
orientation will be denoted by $\vec N^{\vec a}_\tau$.

\subsection{The Jacobi operator}

Let $\Sigma$ be a constant mean curvature surface. Any surface
which is close to $\Sigma$ may be represented as a normal graph
over $\Sigma$
\[
\Sigma_w = \{x + w(x) \, \vec N(x)\quad : \quad  x \in \Sigma\},
\]
where $\vec N$ is the unit normal vector field compatible with the
orientation of $\Sigma$ and $w$ is a (small) scalar function. The
mean curvature of $\Sigma_w$ is denoted by $H(w)$. By definition,
the Jacobi operator about $\Sigma$ is the differential of the
mapping $w \longrightarrow 2 \, H(w)$ at $w=0$. It is given by
\[
\calL_\Sigma := \Delta_\Sigma + |A_\Sigma|^2.
\]
A solution $w$ of the homogeneous problem $\calL_\Sigma w = 0$ is
called a Jacobi field.

\medskip

We denote by $\calL_{D_\tau}$ the Jacobi operator associated to
the Delaunay surface $D_\tau$. In terms of the isothermal
parametrization given in the previous paragraph, it is given by
\[
\calL_{D_\tau} = \frac{4}{\tau^2e^{2\sigma}}\left(\del_{s}^2 +
\del_{ \theta}^2 + \tau^2 \, \cosh (2\sigma) \right).
\]
For the sake of simplicity, we shall now assume that $\tau \neq
1$, namely that $D_\tau$ is not a cylinder. There is no loss of
generality in doing so since our construction, which is based on
the end-to-end construction, does not work for surfaces which have
ends asymptotic to cylinders. Some Jacobi fields are easy to
describe since they correspond to explicit geometric deformations
of the Delaunay surfaces \cite{Maz-Pac}. We briefly describe these
now since they will play a key role in the subsequent analysis.

\medskip

The Jacobi fields corresponding to an infinitesimal translations
of $D_\tau$ will be denoted by $\Phi_\tau^{T, \vec e}$, where
$|\vec e|=1$. They are obtained by projecting the constant vector
field $\vec e$ on the normal vector field $\vec N_\tau$ on
$D_\tau$.
\[
\Phi_\tau^{T,\vec e} : = \vec e\cdot \vec N_\tau
\]
It is geometrically obvious that there are $3$ linearly
independent such Jacobi fields (this is where we use the fact that
$\tau \neq 1$ and hence $D_\tau$ is not a cylinder) which only
depend only on $s$ and are periodic, hence which are bounded as $s
\to \pm \infty$.

\medskip

The two Jacobi fields corresponding to infinitesimal rotations of
the axis of $D_\tau$ will be denoted by $\Phi_\tau^{R,\vec e}$,
where $|\vec e|=1$. They obtained  by projecting the Killing
vector fields
\[
(x_1, x_2, x_3) \longrightarrow  (\vec x \cdot \vec e \, ' ) \,
\vec e \,'' - (\vec x \cdot \vec e \,'') \, \vec e \,'
\]
where $(\vec e, \vec e \, ', \vec e \, '')$ is a direct
orthonormal basis, on the vector field $\vec N_\tau$.
\[
\Phi_\tau^{R, \vec e} : = (\vec x \cdot \vec e \,' ) \, \vec e \,
'' \cdot \vec N_\tau- (\vec x \cdot \vec e \, '') \, \vec e \, '
\cdot \vec N_\tau.
\]
It is geometrically obvious that there are $2$ linearly
independent such Jacobi fields which only depend only on $s$ and
grow linearly in $s$.

\medskip

So far all the Jacobi fields we have defined can be explicitly
computed in terms of the function $\sigma$ and its derivatives
\cite{Maz-Pac}, even though we will not need these expressions.
There is a last Jacobi field, whose geometric meaning is obvious
but whose analytical expression is more intricate, which will be
denoted by $\Phi^{D}_\tau$ and which corresponds to the one
parameter family $D_\tau$ obtained by varying the Delaunay
parameter $\tau$. Since $D_\tau$ are surfaces of revolution, this
Jacobi field depends only on $s$. The fact that $\del_\tau T_\tau
\neq 0$ when $\tau \neq 1$ implies that this Jacobi field is
linearly growing in $s$. Observe that, there exists $p_\tau \in
{\R}-\{0\}$ such that
\begin{equation}
\Phi^{D}_\tau (\cdot + 2 \, s_\tau) = \Phi^{D}_\tau + p_\tau \,
\Phi^{T}_\tau \label{shift} \end{equation} This follows from the
fact that $s \longrightarrow \Phi^{D}_\tau (s + 2 \, s_\tau) -
\Phi^{D}_\tau (s)$ is a bounded Jacobi field which only depends on
$s$, hence is proportional to $\Phi^T_\tau$. The constant $p_\tau$
is not equal to $0$ since $\Phi^{D}_\tau$ is linearly growing.

\medskip

The Jacobi operator ${\cal L}_{D_\tau}$ being invariant with
respect to rotations about the Delaunay axis, we can perform the
eigenfunction decomposition of any function $(s, \theta)
\longrightarrow w(s, \theta)$ in the $\theta$ variable and the
analysis of ${\cal L}_{D_\tau}$ reduces to the analysis of the
sequence of operators
\[
L_{\tau,j}:= \frac{4}{\tau^2 \, e^{2\sigma}} \, \left( \del_{s}^2 +
\tau^2\,\cosh(2\sigma)-j^2 \right)
\]
for $j \in {\mathbb N}$. The potential in $L_{\tau,j}$ being
periodic of period $s_\tau$ (observe that $\cosh (2\sigma)$ is
$s_\tau$ periodic since $\sigma$ is $2 \, s_\tau$ periodic and
odd), it follows from Bloch wave theory \cite{Maz-Pac-Pol-2} that
the following alternative holds~:
\begin{itemize}
\item[(i)] Either the homogeneous problem $L_{\tau,j} w=0$ has two
independent solutions $w^{\pm}$ (depending on $\tau$ and $j$) such
that
\[
w^{\pm}(s+s_\tau) = e^{\pm \zeta_{\tau,j} \, s_\tau} \,
w^{\pm}(s).
\]
for some complex number $\zeta_{\tau ,j}$, with $\Re \,
\zeta_{\tau,j} \geq 0$.

\item[(ii)] Or the homogeneous problem $L_{\tau,j} w=0$ has one
periodic solution and one linearly growing solution. In which
case, we set $\zeta_{\tau ,j} : = 0$.
\end{itemize}

For each $j$, we define the indicial roots associated to the
operator $L_{\tau, j}$ to be the real numbers $\pm
\gamma_{\tau,j}$ where
\[
\gamma_{\tau,j} : = \mbox{Re\,} \, \zeta_{\tau,j} \geq 0.
\]
It is proven in \cite{Maz-Pac-Pol-2} that~:
\begin{proposition}
The indicial roots of $L_\tau$ satisfy the following properties~:
\begin{itemize}
\item[(i)] For any $\tau \in (-\infty,0) \cup (0,1]$,
$\gamma_{\tau,0} =\gamma_{\tau,1} =0$.
\item[(ii)] There exists $\tau_* <0$ such that, for all $j \geq 2$ and
$\tau \in (\tau_*,0) \cup (0,1]$,
$\gamma_{\tau,j} > 0$.
\end{itemize}
\label{pr:3.1}
\end{proposition}
The first property is a consequence of the fact that the Jacobi
fields $\Phi^D_\tau$, $\Phi^{T, \vec e}_\tau$ and $\Phi^{R,\vec
e}_\tau$ are either bounded or linearly growing.

\section{Moduli space theory}

We now briefly describe the moduli space theory for $k$-ended
complete noncompact constant mean curvature surfaces as developed
in \cite{Kus-Maz-Pol} and extended in \cite{Maz-Pac-Pol-2}. We
define $\calM_{g,k}^{\tau_*}$ to be the set of all complete,
noncompact constant mean curvature surfaces which have genus $g$
and $k$ ends asymptotic to Delaunay surfaces whose Delaunay
parameter belongs to $(\tau_*,0)\cup(0,1]$. Observe that we do not
mod out by the group of rigid motions.

\medskip

We can decompose a surface $\Sigma \in \calM_{g,k}^{\tau_*}$ into
overlapping connected pieces~: A compact component $K$ and the
ends $E_\ell$, for $\ell=1, \ldots, k$ and we can require that
each $K \cap E_\ell$ is homeomorphic to an annulus $[0,1]\times
S^1$. For each $\ell$, we choose standard isothermal coordinates
$(s,\theta)$ for the model Delaunay end $D_{\tau_\ell}$ so that
the end $E_\ell$ is parametrized by
\begin{equation}
Y_\ell : =  X^{\vec a_\ell}_{\tau_\ell} + w_\ell \, \vec
N_{\tau_\ell}^{\vec a_\ell} + \vec b_\ell , \label{eq:padam}
\end{equation}
for $(s, \theta) \in [0 ,+\infty) \times S^1$. Since we have
assumed that the end $E_\ell$ is asymptotic to $D^{\vec
a_\ell}_{\tau_\ell}+ \vec b_\ell$, this means that the function
$w_\ell$ is exponentially decreasing. To be more specific, we
 need the~:
\begin{definition}
Given $r \in \N$, $\alpha \in (0,1)$ and $\mu \in {\R}$, the space
${\mathcal E}^{r,\alpha}_\mu([0, +\infty) \times S^1)$ is the
space of functions $v \in{\calC}^{r,\alpha}_{loc} ([0, +\infty)
\times S^1)$ for which
\[
\| v \|_{\calE^{r,\alpha}_\mu ([0, +\infty) \times S^1)} : =
\sup_{s \geq 0} e^{-\mu \, s} \, \left| v
\right|_{\calC^{r,\alpha}([s,s+1] \times  S^1)}
\]
is finite. \label{de:1}
\end{definition}
Granted this definition, it is known that
\begin{equation}
w_\ell \in \calE^{2,\alpha}_{-
\gamma_{\tau_\ell,2}}([0,\infty)\times S^1). \label{eq:6.22}
\end{equation}
In other words the rate of decay of the function $w_\ell$ is
dictated by the indicial root $\gamma_{\tau_\ell ,2}$. We refer to
\cite{Maz-Pac} for a proof of this fact. The moduli space theory
is based on the~:
\begin{definition}
For $r \in \N$, $\alpha \in (0,1)$ and $\mu \in \R$, let
${\mathcal D}^{r,\alpha}_\mu(\Sigma)$ be the space of functions $v
\in{\calC}^{r,\alpha}(\Sigma)$ for which
\[
\| v \|_{\calD^{r,\alpha}_\mu} : = \|\left. v
\right|_K\|_{\calC^{r, \alpha}} + \sum_{\ell = 1}^k \|\left. v
\circ Y_\ell \right|_{E_\ell} \|_{\calE^{ r,\alpha}_\mu}
\]
is finite. \label{de:2}
\end{definition}
We can now give the precise definition of a nondegenerate constant
mean curvature surface.
\begin{definition}
The surface $\Sigma \in \calM_{g,k}^{\tau_*}$ is nondegenerate if
\[
\calL_{\Sigma} : {\mathcal D}^{2,\alpha}_\mu(\Sigma)
\longrightarrow {\mathcal D}^{0,\alpha}_\mu(\Sigma)
\]
is injective for all $\mu < 0$.
\label{de:7.2}
\end{definition}

Following the analysis of the Jacobi fields we have done in \S 2.2
and using the parameterization (\ref{eq:padam}) together with
(\ref{eq:6.22}), it is easy to see that, on each end $E_\ell$ of
$\Sigma$, there exists $5$ (globally defined) independent Jacobi
fields $\Phi_{E_\ell}^{T,\vec e}$ and $\Phi_{E_\ell}^{R,\vec e}$
which satisfy
\begin{equation}
\begin{array}{rlllll}
\Phi^{T,\vec e }_{E_\ell} \circ Y_\ell - \Phi^{T,\vec
e}_{\tau_\ell} \in \calE^{2,\alpha}_{-\gamma_{\tau_\ell,2}} ([0 ,
+\infty) \times S^1)), \\[3mm]
\Phi^{R,\vec e }_{E_\ell} \circ Y_\ell - \Phi^{R,\vec
e}_{\tau_\ell} \in \calE^{2,\alpha}_{-\gamma_{\tau_\ell,2}} ([0 ,
+\infty) \times S^1)), \label{eq:est}
\end{array}
\end{equation}
where $|\vec e|=1$.

\medskip

The existence of a Jacobi field $\Phi^{D}_{E_\ell}$ (only defined
on $E_\ell$) which is asymptotic to $\Phi^{D}_\tau$ is not a
trivial fact. This follows from a perturbation argument
\cite{Maz-Pac} and, in general, this Jacobi field is only defined
on $E_\ell$ away from a compact set in $\Sigma$ and is not
globally defined. This motivates the~:
\begin{definition}
The end $E_\ell$ of $\Sigma$ is said to be regular if there exists
a globally defined Jacobi field $\Phi^{D}_{E_\ell}$ satisfying
\begin{equation}
\Phi^{D}_{E_\ell} \circ Y_\ell - \Phi^{D}_{\tau_\ell} \in
\calE^{2,\alpha}_{\mu} ([0 , +\infty) \times S^1)),
\end{equation}
for all $\mu \in (-\gamma_{\tau_\ell,2}, 0)$.
\end{definition}
The fact that such a globally defined Jacobi field exists is
usually a consequence of the existence of a one parameter family
of constant mean curvature surfaces $\Sigma (\varepsilon)$, for
$\varepsilon \in (-\varepsilon_0,\varepsilon_0)$, which have $k$
ends, are close to $\Sigma$ (in a suitable sense), satisfy
$\Sigma_0 =\Sigma$ and whose $\ell$-th end $E_\ell (\varepsilon)$
is asymptotic to a Delaunay surface of parameter $\tau_\ell +
\varepsilon$.

\medskip

As in \cite{Kus-Maz-Pol}, we define the $6k$-dimensional
deficiency space
\[
{\mathcal W}_\Sigma : =  \oplus_{\ell=1}^k \mbox{Span} \left\{
\chi_{E_\ell} \, \Phi^{D}_{E_\ell}  ,  \chi_{E_\ell} \, \Phi^{T,
\vec e }_{E_\ell} , \chi_{E_\ell} \, \Phi^{R,\vec e }_{E_\ell},
\quad : \quad |\vec e| = 1 \right\},
\]
where $\chi_{E_\ell}$ is a cutoff function equal to $0$ on $
\Sigma -E_\ell$ and equal to $1$ on $Y_{\ell}([1, \infty)\times
S^1)$. The following Proposition is the key result for the study
of the structure of ${\calM}^{\tau_*}_{g,k}$.
\begin{proposition}
\cite{Kus-Maz-Pol} Assume that $\Sigma \in {\calM}^{\tau_*}_{g,k}$
is nondegenerate and fix $\mu \in (-\inf_\ell
\gamma_{\tau_\ell,2},0)$. Then the mapping
\begin{equation}
\calL_\Sigma : \calD^{2,\alpha}_\mu(\Sigma) \oplus \calW_\Sigma
\longrightarrow \calD^{0,\alpha}_\mu(\Sigma)
\label{eq:7.1a}
\end{equation}
is surjective and has a kernel of dimension $3k$. Moreover, there
exists a $3k$-dimensional subspace $\calN_\Sigma \subset
\calW_\Sigma$ such that
\[
\mbox{Ker} \, \calL_\Sigma \subset \calD^{2,\alpha}_\mu(\Sigma)
\oplus \calN_\Sigma
\]
Finally, given any $3k$-dimensional subspace $\calK_\Sigma \subset
\calW_\Sigma$ such that $ \calK_\Sigma \oplus \calN_\Sigma =
\calW_\Sigma$ the mapping
\begin{equation}
\calL_\Sigma :\calD^{2,\alpha}_\mu(\Sigma) \oplus \calK_\Sigma
\longrightarrow \calD^{0,\alpha}_\mu(\Sigma)
\label{eq:7.1}
\end{equation}
is an isomorphism.
\label{pr:7.1}
\end{proposition}
It follows from this result that ${\calM}^{\tau_*}_{g,k}$ is
locally a $3k$-dimensional smooth manifold near any nondegenerate
element \cite{Kus-Maz-Pol} (observe that we have not taken the
quotient by the group of rigid motions of ${\R}^3$).

\section{Building blocks}

We describe two families of complete noncompact constant mean
curvature surfaces which will be used in the construction. The
members of the first family are $3$-ended surfaces while the
members of the second family are $k$-ended surfaces. We give a
fairly precise description of the elements of each family and
explain how these families can be obtained using already known
constructions of complete noncompact constant mean curvature
surfaces. In this paper we do not give a proof of the existence of
these families but rather to rely on their existence. We hope that
the reader will either be convinced by the explanations below or
take the existence of these families for granted.

\medskip

We start by recalling the well known balancing formula
\cite{Kor-Kus-Mee}. Given a constant mean curvature surfaces
$\Sigma\subset \R^{3}$ with finitely many ends $E_\ell$, for $\ell
=1, \ldots, k$, which are asymptotic to Delaunay surfaces ${\cal
D}_{\tau_\ell}^{\vec a_\ell} + \vec b_\ell$, the balancing formula
reads~:
\begin{equation}
\sum_{\ell =1}^k \tau_\ell  \, |\tau_\ell | \, \vec{a}_\ell =0
\label{BF}
\end{equation}
where $\vec{a}_\ell$ is the direction of the axis of $E_\ell$,
which is normalized by $|\vec{a}_\ell|=1$ and points toward the
end of $E_\ell$.

\medskip

We fix $(\vec e_1, \vec e_2, \vec e_3)$ a direct orthonormal basis
of ${\R}^3$.

\subsection{Type-1 surfaces}

The members of the first family are denoted by $\Sigma_{\tau
,\alpha}$, where $\tau$ and $\alpha$ are parameters. These
surfaces are assumed to enjoy the following properties~:

\begin{itemize}
\item[(i)] Each $\Sigma_{\tau , \alpha}$ is a complete noncompact
constant mean curvature surface with $3$ ends which are denoted by
$E^{-1}_{\tau ,\alpha}, E^{0}_{\tau ,\alpha}$ and $E^{1}_{\tau
,\alpha}$.

\item[(ii)] The surface $\Sigma_{\tau , \alpha}$ is invariant
under the action of the group
\[
G : = \{ I, S_1, S_3\}
\]
where $S_i$ is the symmetry with respect to the plane $x_i=0$.

\item[(iii)] Each $\Sigma_{\tau , \alpha }$ is nondegenerate and
the parameters $(\tau , \alpha)$ are local parameters on the
moduli space of constant mean curvature surfaces with $3$ ends,
which are invariant under the action of the group $G$.

\item[(iv)] The end $E^0_{\tau , \alpha}$ is asymptotic to a
Delaunay surface of parameter $\tau$ and axis the $x_2$-axis. The
vector $- \vec{e}_2$ is directed toward the end of $E^0_{\tau ,
\alpha}$. In particular, there exists a smooth function $\tau
\longrightarrow d_{\tau , \alpha}^0$ such that $E^0_{\tau ,
\alpha}$ is a graph (for an exponentially decaying function) over
the Delaunay surface $D_{\tau}^{{\vec e}_2} - d_{\tau , \alpha}^0
\, \vec e_2$.

\item[(v)] The end $E^1_{\tau , \alpha}$ is asymptotic to the
Delaunay surface of parameter $\bar \tau$ and axis passing through
the origin and of direction
\[
\vec{a}_\alpha := - \sin \alpha \, \vec{e}_1 - \cos \alpha \,
\vec{e_2} .
\]
The vector $\vec{a}_\alpha$ is directed toward the end of
$E^1_{\tau , \alpha}$. In particular, there exists a smooth
function $\tau \longrightarrow d_{\tau , \alpha}^1$ such that
$E^1_{\tau , \alpha}$ is a graph (for an exponentially decaying
function) over the Delaunay surface $D_{\bar \tau}^{{\vec
a}_\alpha} + d_{\tau , \alpha}^1 \, \vec{a}_\alpha$.

\item[(vi)] The ends of $\Sigma_{\tau , \alpha}$ are regular.
\end{itemize}

Observe that the image of $E^1_{\tau , \alpha}$ by $S_1$ is
$E^{-1}_{\tau , \alpha}$ and that $E^0_{\tau , \alpha}$ remains
globally fixed under the action of $S_1$. Also each end remains
globally fixed under the action of $S_3$. Applying the balancing
formula (\ref{BF}), we conclude that the Delaunay parameters $\bar
\tau$ and $\tau$ are related by the formula
\begin{equation}
\tau \, |\tau| + 2 \, \cos \alpha \, \bar \tau  \, |\bar \tau |=
0.
\label{Bla}
\end{equation}
In particular, if $\alpha \in (0, \pi /2)$, the signs of $\tau$
and $\bar \tau$ are different and this implies that the surface
$\Sigma_{\tau , \alpha}$ has always an end which is not embedded
(asymptotic to a nodoid) in this case.

\medskip

Observe that (iv) implies that the end $E^0_{\tau , \alpha}$ can
be parameterized by
\begin{equation}
X^0_{\tau, \alpha} (s, \theta) : = X_\tau^{\vec e_2} (s, \theta) -
d_{\tau , \alpha}^0 \, \vec e_2 + w^0_{\tau, \alpha} (s, \theta)
\, \vec N^{\vec e_2}_\tau (s, \theta) \label{eq:dldl0}
\end{equation}
with $(s, \theta) \in [0, \infty) \times S^1$, for some function
$w^0_{\tau, \alpha} \in {\cal E}^{2, \alpha}_{-\gamma_{\tau,2}}
([0, +\infty) \times S^1)$ (In general the function $w^0_{\tau,
\alpha}$ is only defined on $[c, + \infty) \times S^1$ for some $c
>0$ large enough. However increasing the value of $d_{\tau,
\alpha}$ by a $2 \,m \, T_\tau$ for some $m \in {\mathbb N}$, if
this is necessary, we can assume that the function $w^0_{\tau
,\alpha}$ is defined on $[0, + \infty) \times S^1$).

\medskip

Similarly (v) implies that the end $E^1_{\tau , \alpha}$ can be
parameterized by
\begin{equation}
X^1_{\tau, \alpha} (s, \theta) : = X_{\bar \tau}^{\vec a_\alpha}
(s, \theta) + d_{\tau , \alpha}^1 \, \vec a_\alpha + w^1_{\tau,
\alpha} (s, \theta) \, \vec N^{\vec a_\alpha }_{\bar \tau} (s,
\theta) \label{eq:dldl1}
\end{equation}
with $(s, \theta) \in [0, + \infty) \times S^1$, for some function
$w^1_{\alpha, \tau} \in {\cal E}^{2, \alpha}_{-\gamma_{\bar \tau ,
2}} ([0, + \infty) \times S^1)$.

\begin{definition}
Given $s_0, s_1 >0$, we define the compact surface with $3$
boundaries
\[
\begin{array}{rllll}
\Sigma_{\tau , \alpha} (s_0, s_1) : = \Sigma_{\tau, \alpha} - (
X_{\tau , \alpha }^0 ((s_0, +\infty) \times S^1) & \cup & X_{\tau,
\alpha}^1 ((s_1, + \infty) \times S^1) \\[3mm]
& \cup & S_1 \, X_{\tau, \alpha}^1 ((s_1, + \infty) \times S^1) )
\end{array}
\]
\label{de:trunc1}
\end{definition}

In the case where the surfaces are Alexandrov embedded the
surfaces described above have been classified in \cite{Kar-Kus}.
However, it does not follow from this description that the
surfaces are nondegenerate and have regular ends. This is the
reason why we give now two examples of construction of such a
family which rely on connected sum constructions and for which it
is possible to check that the surfaces constructed are both
nondegenerate and have regular ends~:

\medskip

\noindent {\bf Example 1} A first family can be obtained by gluing
on the unit sphere $S^2 \subset \R^{3}$, three half Delaunay
surfaces of parameters $\bar \tau$, $\tau$ and $\bar \tau$
respectively at the points of coordinates
\[
(-\sin \alpha, - \cos \alpha , 0,),  \qquad (0,-1, 0,) \quad
\mbox{and} \quad (\sin \alpha, - \cos \alpha , 0)
\]
respectively, using a modified version of the connected sum result
of \cite{Maz-Pac-Pol-1},\cite{Maz-Pac-Pol-2} and \cite{Jle}. The
construction works if one imposes the surfaces to be invariant
under the action of the group $G$. Given the symmetries of the
surfaces constructed, there remains only two degrees of freedom
which are : The Delaunay parameter $\tau$ and the angle $\alpha$
between the ends. The construction works for any $\alpha \in (0,
\pi/2) \cup (\pi/2, \pi)$ and any $\tau \neq 0$ close enough to
$0$. The fact that the ends are regular follows from the
construction itself since $\tau$ can be used to parameterize this
family of surfaces and differentiation with respect to this
parameter yields a Jacobi field whose asymptotic along any end has
a nontrivial component on $\chi_{E_\ell} \, \Phi^D_{E_\ell}$, for
$\ell=0, \pm 1$.

\medskip

\noindent {\bf Example 2} A second family can be obtained by
gluing on a Delaunay surface of parameter $\bar \tau$ and axis
$x_1$ which is translated so that it is invariant under the action
of the symmetry $S_1$ (namely either $D^{\vec e_1}_{\bar \tau}$ or
$D^{\vec e_1}_{\bar \tau} + T_{\bar \tau} \, \vec e_1$), a half
Delaunay surface of axis $x_2$ and small Delaunay parameter
$\tau$. Again, the construction works if one imposes the surfaces
to be invariant under the action of the group $G$. Given the
symmetries of the surfaces constructed, there remains only two
degrees of freedom which are : The Delaunay parameters $\bar \tau$
and $\tau$. The construction works for any small value of the
parameter $\tau \neq 0$ \cite{Maz-Pac-Pol-1}, \cite{Maz-Pac-Pol-2}
and \cite{Jle} and provides a surface with an angle $\alpha $
close, but not equal, to $\pi/2$ which is determined by the
equation $\tau \, |\tau| + 2\,\cos \alpha \, \bar \tau \, |\bar
\tau |= 0$. This shows that $(\alpha, \tau)$ are local parameters
on the corresponding moduli space and, as in the previous example,
the ends of the surfaces are regular.

\medskip

In both cases, the surfaces are seen to be nondegenerate, when
$\tau$ is close enough to $0$, using the strategy developed in
\cite{Maz-Pac}.

\subsection{Type-$2$ surfaces}

We fix $k \geq 3$. The members of the second family are denoted by
$\bar \Sigma_\tau$, where $\tau$ is a parameter. These surfaces
are assumed to enjoy the following properties~:

\begin{itemize}
\item[(i)]  Each $\bar \Sigma_\tau$ is a complete noncompact
constant mean curvature surface with $k$ ends which are denoted by
$ \bar E^0_\tau, \ldots,  \bar E^{k-1}_\tau $.

\item[(ii)] The surface in invariant under the action of the group
\[
G_k : = \{ R_{2\pi j/k} \quad : \quad j \in {\mathbb Z} \}
\]
where $R_{\theta}$ is the rotation of angle $\theta$ in the $x_1, x_2$
plane.

\item[(iii)] Each $\bar \Sigma_\tau$ is nondegenerate and the
parameter $\tau$ is local parameter on the moduli space of
constant mean curvature surfaces with $k$ ends, which are
invariant under the action of the group $G_k$.

\item[(iv)] The end $\bar E^0_\tau$ is asymptotic to a Delaunay
surface of parameter $\tau$ and axis the $x_2$-axis. The vector
$\vec{e}_2$ being directed toward the end of $\bar E^0_\tau$. In
particular, there exists a smooth function $\tau \longrightarrow
\bar d^0_\tau$ such that $\bar E^0_\tau$ is a graph (for an
exponentially decaying function) over the Delaunay surface
$D_{\tau}^{{\vec e}_2} + \bar d^0_{\tau} \, \vec e_2$.

\item[(v)] The ends of $\bar \Sigma_\tau$ are regular.
\end{itemize}

Observe that, for $\ell=1, \ldots, k-1$ the image of $E^0_\tau$ by
$R_{2\pi \ell/k}$ is the end $\bar E^\ell_\tau$. Hence the angle
between two consecutive ends is given by  $2\, \pi/k$ and, to
check that the ends of $\bar \Sigma_\tau$ are regular it is enough
to check that $\bar E_\tau^0$ is regular.

\medskip

As in the case of Type-1 surfaces, (iv) implies that the end $\bar
E^0_\tau$ can be parameterized by
\begin{equation}
\bar X_{\tau}^{0} (s, \theta) : = X_\tau^{\vec e_2} (s, \theta) +
\bar d^0_\tau \, \vec e_2 + \bar w^0_\tau (s, \theta) \, \vec
N^{\vec e_2}_\tau (s, \theta) \label{eq:ldld}
\end{equation}
with $(s, \theta) \in (0, +\infty) \times S^1$, for some function
$\bar w^0_\tau \in {\cal E}^{2, \alpha}_{-\gamma_{\tau , 2}} (0,
+\infty)$.

\begin{definition}
Given $s_0 >0$, we define the compact surface with $k$ boundaries
\[
\bar \Sigma_{\tau} (s_0) : = \bar \Sigma_\tau -
\cup_{\ell=0}^{k-1} R_{2\pi\ell/k} \, \bar X_{\tau}^0 ((s_0,
+\infty) \times S^1)
\]
\label{de:trunc2}
\end{definition}
We now give two examples of such a family.

\medskip

\noindent {\bf Example 1} A first family can be obtained by gluing
on the unit sphere $S^2\subset\R^{3}$, $k$ copies of a half
Delaunay surface with small Delaunay parameter $\tau \neq 0$ in
such a way that the surface remains invariant under the action of
$G_k$. Again this is a byproduct of the end addition result proved
in \cite{Maz-Pac-Pol-1}, \cite{Maz-Pac-Pol-2} or this is also a
byproduct of the result of N. Kapouleas in \cite{Kap-1}. These
surfaces have also been constructed and described by K.
Grosse-Brauckmann \cite{Kar}.

\medskip

\noindent {\bf Example 2} A second family can be obtained by
gluing on a $k$-noid (a minimal surface with $k$ ends of
catenoidal type \cite{Jor-Mee}, \cite{Cos-Ros}) which is invariant
under the action of $G_k$, $k$ copies of a half Delaunay surface
with small Delaunay parameter $\tau \neq 0$ in such a way that the
symmetries are preserved. This construction is the one described
in \cite{Maz-Pac}.

\medskip

In either case, given the symmetries of the surfaces constructed
there remains only one degree of freedom which is $\tau$, the
Delaunay parameter of the ends. Either construction works for any
$\tau \neq 0$ close enough to $0$. The fact that (v) holds follows
at once from the construction itself since $\tau$ can be used to
parameterize this family of surfaces and differentiation with
respect to this parameter yields a Jacobi field whose asymptotic
has a nontrivial component on $\chi_{\bar E_\tau^0} \,
\Phi^D_{\bar E^0_\tau}$. The fact that the surfaces constructed
are nondegenerate follows from \cite{Maz-Pac}.

\subsection{Jacobi fields}

We give a precise description of the Jacobi fields on both
$\Sigma_{\tau , \alpha}$ and on $\bar \Sigma_\tau$. This
description yields a description of the spaces ${\cal
K}_{\Sigma_{\tau , \alpha}}$ and ${\cal K}_{\bar \Sigma_\tau}$
which have been introduced in Proposition~\ref{pr:7.1}.

\medskip

We start with the analysis  of  the Jacobi fields on $\bar
\Sigma_\tau$ since this is the simplest. Since the surface $ \bar
\Sigma_\tau$ is assumed to be nondegenerate, the deficiency space
${\mathcal D}_{\bar \Sigma_\tau }$ is $6k$-dimensional. However,
since we are working in the space of surfaces which are invariant
under the action of the group $G_k$ and this reduces the dimension
of the corresponding moduli space to $1$ and the deficiency space
is now spanned by the $2$ functions
\[
\bar \psi_\tau^{T} : = \sum_{\ell =0}^{k-1} \chi_{\bar
E^\ell_{\tau}} \, \Phi_{\bar E^\ell_\tau }^{T, \vec a_\ell} \qquad
\mbox{and} \qquad \bar \psi_\tau^{D} : = \sum_{\ell =0}^{k-1}
\chi_{\bar E^\ell_{\tau}} \, \Phi_{\bar E^\ell_\tau }^{D}
\]
where $\vec a_\ell = R_{2\pi\ell/k} \, \vec e_2$ is the direction
of the end $\bar E^\ell_\tau$. Observe that the symmetries of
$\Sigma_\tau$ imply that
\[
\Phi_{\bar E^\ell_\tau }^{T, \vec a_\ell} = \Phi_{\bar E^0_\tau
}^{T, \vec e_2} \circ  (R_{2\pi\ell/k})^{-1}
\]
Since the end $\bar E^0_\tau$ is assumed to be regular, there
exists a globally Jacobi field (which is invariant under the
action of $G_k$) whose asymptotic on $\bar E^0_\tau$ has a
nontrivial component on $\Phi^D_{\bar E^0_\tau}$. In fact this
Jacobi field is obtained by moving the parameter $\tau$.
Multiplying this Jacobi field by a suitable constant, we can
assume that it is asymptotic to $\bar \psi^D_\tau + \bar c \, \bar
\psi^T_\tau$ on each $\bar E^\ell_{\tau}$, were the constant $\bar
c$ depends on $\tau$. This implies that the space ${\mathcal
K}_{\bar \Sigma_\tau }$ can be chosen to be
\[
{\mathcal K}_{ \bar \Sigma_\tau } = \mbox{Span} \{ \bar
\psi^T_\tau \}
\]

We now analyze the Jacobi fields on $\Sigma_{\tau, \alpha}$. By
assumption, $\Sigma_{\tau , \alpha}$ is nondegenerate and has $3$
ends, therefore the deficiency space ${\mathcal D}_{\Sigma_{\tau ,
\alpha}}$ is $18$-dimensional. Now, recall that we are working in
the space of surfaces which are invariant under the action of the
group $G$ and this reduces the dimension of the corresponding
moduli space to $3$ and the deficiency space is spanned by the $6$
functions we now describe~:
\[
\begin{array}{lllllll}
\psi_{\tau, E^0}^{T} : = \chi_{E^0_{\tau, \alpha}} \,
\Phi_{E^0_{\tau, \alpha}}^{T, \vec e_2},

&\quad  \psi_{\tau , E^0}^{D} : =  \chi_{E^0_{\tau, \alpha}} \,
\Phi_{E^0_{\tau, \alpha}}^D ,\\[3mm]

\psi_{\tau, E^1}^{D} : = \chi_{E^1_{\tau, \alpha}} \,
\Phi_{E^1_{\tau, \alpha}}^{D} +  (\chi_{E^1_{\tau, \alpha}}
\Phi_{E^1_{\tau, \alpha}}^{D}) \circ S_1,

& \quad \psi_{\tau, E^1}^{T, \vec a} : = \chi_{E^1_{\tau, \alpha}}
\, \Phi_{E^1_{\tau, \alpha}}^{T,\vec a_\alpha } +
(\chi_{E^1_{\tau, \alpha}} \Phi_{E^1_{\tau, \alpha}}^{T, \vec
a_\alpha})\circ S_1, \\[3mm]

\psi_{\tau,  E^1}^{T, \vec a^\perp} : =\chi_{E^1_{\tau, \alpha}}
\, \Phi_{E^1_{\tau, \alpha}}^{T,\vec a_\alpha^\perp } +
(\chi_{E^1_{\tau, \alpha}} \Phi_{E^1_{\tau, \alpha}}^{T, \vec
a_\alpha^\perp}) \circ S_1,

& \quad  \psi_{\tau, E^1}^{R} : = \chi_{E^1_{\tau, \alpha}} \,
\Phi_{E^1_{\tau, \alpha}}^{R,\vec e_3 } + (\chi_{E^1_{\tau,
\alpha}} \Phi_{E^1_{\tau, \alpha}}^{R, \vec e_3}) \circ S_1,
\end{array}
\]
where $\vec a : = \vec a_\alpha$ and $\vec a^\perp : = \vec
a_\alpha^\perp : = \cos \alpha \, \vec e_1+ \sin \alpha \, \vec
e_2$. Even though these functions do depend on $\alpha$, we have
not indicated this.

\medskip

We now describe the Jacobi fields which are globally defined on
$\Sigma_{\tau, \alpha}$ since they are obtained by moving the two
parameters $\alpha$, $\tau$ and also by translating this surface
in the $\vec e_2$ direction. These Jacobi fields span the
nullspace ${\cal N}_{\Sigma_{\tau, \alpha}}$.
\begin{itemize}
\item[(1)] Changing the $\tau$ parameter (keeping $\alpha$ fixed)
changes $\bar \tau$. Therefore, this yields a Jacobi field which
(up to a multiplicative constant) is asymptotic to
\[
\psi_{\tau, E^1}^{D} + c_1 \, \psi_{\tau, E^1}^{T, \vec a}
\]
on $E^1_{\tau, \alpha} \cup E^{-1}_{\tau, \alpha}$ and which is
asymptotic to $c_2 \, \psi_{\tau, E^0}^{D} + c_3 \, \psi_{\tau,
E^0}^{T}$ on $E^0_{\tau, \alpha}$.

\item[(2)] Changing the $\alpha$ parameter (keeping $\tau$ fixed),
yields a Jacobi field which (up to a multiplicative constant) is
asymptotic to
\[
\psi_{\tau ,E^1}^{R} + c_4 \, \psi_{\tau, E^1}^{T,\vec a} + c_5 \,
\psi_{\tau, E^1}^{D} + c_6 \, \psi_{\tau , E^1 }^{T, \vec a^\perp}
\]
on $E^1_{\tau, \alpha} \cup E^{-1}_{\tau, \alpha}$ and which is
asymptotic to $c_7 \, \psi_{\tau, E^0}^{T}$ on $E^0_{\tau,
\alpha}$.

\item[(2)] Translation of $\Sigma_{\tau, \alpha}$ in the $\vec
e_2$ direction (keeping $\tau$ and $\alpha$ fixed) yields a Jacobi
field which is asymptotic to
\[
\psi^{T}_{\tau, E^0}
\]
on $E^0_{\tau, \alpha }$ and which is asymptotic to $ c_8 \,
\psi_{\tau, E^1}^{T, \vec a} + c_9 \, \psi_{\tau, E^1}^{T,\vec
a^\perp}$ on $E^1_{\tau, \alpha}$.
\end{itemize}
Here the constants $c_1, \ldots, c_9$  depend on $\tau$ and
$\alpha$. Recall that the space ${\mathcal K}_{\Sigma_{\tau,
\alpha}}$ is a $3$-dimensional subspace of the deficiency space
${\mathcal D}_{\Sigma_{\tau, \alpha }}$ chosen so that
\[
{\mathcal D}_{\Sigma_{\tau, \alpha }} = {\cal K}_{\Sigma_{\tau,
\alpha }} \oplus {\cal N}_{\Sigma_{\tau, \alpha }}.
\]
It follows from the above description of the elements of
${\mathcal N}_{\Sigma_{\tau, \alpha }}$ that we can choose
\[
{\cal K}_{\Sigma_{\tau, \alpha }} = \mbox{Span} \{ \psi_{\tau,
E^0}^{D} + t \, \psi_{ \tau, E^0}^T, \psi_{\tau, E^1}^{T, \vec a}
, \psi_{\tau, E^1}^{T, \vec a^\perp}\}
\]
where $t \in {\mathbb R}$ is a free parameter which will be fixed
later on.

\section{The construction}

We fix $k \geq 3$ and define
\[
\alpha_k : = \frac{\pi}{2} - \frac{\pi}{k}
\]
We assume that, for $\tau$ in some closed (nonempty) interval $I
\subset (\tau_*, 0) \cup (0, 1)$, we are given a family of
surfaces $\Sigma_{\tau , \alpha_k}$ of Type $1$ and a family of
surfaces $\bar \Sigma_\tau$ of Type $2$. For the sake of
simplicity we now drop the dependence on $\alpha_k$ in all the
quantities related to $\Sigma_{\tau , \alpha_k }$ and simply write
$\Sigma_\tau $, $E_\tau^\ell$, $d^0_\tau$, $d^1_\tau$, \ldots The
parameter $\tau$ being chosen in $I$, we recall that $\bar \tau$
is given by
\begin{equation}
\tau \, |\tau| + 2\,\cos \alpha_k \, \bar \tau \, |\bar \tau| = 0.
\label{eq:tta} \end{equation} Given $n, m \in {\mathbb N}$, we set
\[
\delta_{n,\tau} : = d^0_{\tau} + \bar d_{\tau}^0 + 2 \, n \,
T_{\tau}
\]
We agree on the notation
\[
\Sigma_{n, \tau} : = \Sigma_\tau + \delta_{n, \tau} \, \vec e_2
\]
and the ends of this surface are denoted by
\[
E^j_{n, \tau} : = E^j_\tau + \delta_{n, \tau} \, \vec e_2
\]
and are parameterized by
\[
 X^j_{n, \tau} : = X^j_\tau + \delta_{n, \tau} \, \vec e_2 .
\]
Also we define the truncated surface (see
Definition~\ref{de:trunc1})
\[
\Sigma_{n,\tau} (s_0, s_1) := \Sigma_\tau (s_0, s_1) +
\delta_{n,\tau} \, \vec e_2
\]

With these notations in mind, we consider the truncated surface
$\Sigma_{n,\tau} (n\, s_\tau, m\, s_{\bar \tau} )$  together with
the images of this surface by $R_{2\ell \pi/k}$, for $\ell=1,
\ldots, k-1$ and also the truncated surface $\bar \Sigma_\tau (n
\, s_\tau)$ (see Definition~\ref{de:trunc2}). These surfaces with
boundaries are now connected together using appropriate cutoff
functions, to produce a compact surface which is invariant under
the action of $G_k$. More precisely, for each $\ell =0, \ldots,
k-1$~: The end $\bar E^\ell_\tau$ of $\bar \Sigma_\tau$ can be
connected with the image of $E_{n,\tau}^0$ by $R_{2\pi \ell/k}$
since they are graphs over the same Delaunay surface. And,
provided $n$ and $\tau$ are suitably chosen, the image of the end
$E^1_{n,\tau}$ by $R_{2\pi\ell /k}$ can be connected with the
image of $E^{-1}_{n,\tau}$ by $R_{2\pi(\ell+1) /k}$. We now
describe analytically this procedure. Given the fact that the
surface we want to construct should be invariant under the action
of $G_k$ it is enough to describe the~:

\subsection{Connection of $E^{0}_{n,\tau}$ with $\bar E^{0}_\tau$.}

By construction the ends $E^0_{n,\tau}$ and $\bar E^0_\tau $ are
normal graphs over the same Delaunay surface. Given the
parameterizations defined in (\ref{eq:dldl0}) and (\ref{eq:ldld})
we can connect the two pieces together by considering the
parameterization
\[
Y^0_{n,\tau} (s, \theta) = X^{\vec e_2}_\tau \left( s + n \,
s_\tau , \theta \right)  + \bar d_\tau^0 \, \vec e_2  + \tilde
w_\tau (s, \theta) \, \vec N^{\vec e_2}_\tau \left( s + n \,
s_\tau , \theta \right)
\]
for $(s, \theta) \in (- n \, s_\tau,  n \, s_\tau) \times S^1$
where
\[
\tilde w_\tau (s, \theta) : =  \xi (s) \, \bar w^0_\tau \left(s +
n\, s_\tau , \theta \right) + ( 1 - \xi (s) ) \, w^0_\tau \left(
n\, s_\tau - s, \theta \right)
\]
Here $\xi$ is a cutoff function identically equal to $1$ for $s
\leq -1$ and identically equal to $0$ for $s \geq 1$ and which
satisfies
\[
\xi (-s) = 1 - \xi (s) .
\]
We will denote by $A^0_{n,\tau}$ the image of $(-1, 1) \times S^1$
by $Y_\tau^0$. We define
\[
Y^\ell_{n,\tau} = R_{2\pi\ell/k} \, Y^0_{n,\tau}
\]
for $\ell =1, \ldots, k-1$ which describes the connection of $\bar
E^\ell_\tau$ with the image of $E^0_{n,\tau}$ by $R_{2\pi
\ell/k}$.

\subsection{Connection of $E^{1}_{n,\tau} $ with the image of $
E^{-1}_{n,\tau}$ by $R_{2\pi/k}$.}

We define the plane
\[
\Pi_k : =  \{ x \in {\R}^3 \, : \, \tan (2\pi/k) \, x_2 = - x_1 \}
\]
Observe that the image of $E^1_{n,\tau}$ by the symmetry with
respect to $\Pi_k$ is equal to the image of $E^{-1}_{n,\tau}$ by
$R_{2\pi/k}$. By definition, the end $E^1_{n,\tau}$ is a graph
over the Delaunay surface $D^{\vec a}_{\bar \tau} + d^1_\tau \,
\vec a + \delta_{n,\tau} \, \vec e_2$. Therefore the end
$E^1_{n,\tau}$ and its image by the symmetry with respect to the
plane $\Pi_k$ are normal graphs over the same Delaunay surface if
and only if the Delaunay surface $D^{\vec a}_{\bar \tau} +
d^1_\tau \, \vec a + \delta_{\tau,n} \, \vec e_2$ is invariant
under the symmetry with respect to the plane $\Pi_k$. This
condition is translated into the fact that there exists an integer
$m \in {\mathbb N}$ such that
\begin{equation}
\sin (\pi/k) \, \left( d^0_\tau + \bar d^0_\tau + 2 \, n \,
T_{\tau} \right) = d^1_\tau + m \, T_{\bar \tau}
\label{eq:Zcondition}
\end{equation}
If this condition is fulfilled we can connect the end
$E^1_{n,\tau}$ and its image by $R_{2\pi/k}$, using the
parameterization
\[
Z^0_{n,\tau} (s, \theta) =  X^{\vec a }_{\bar \tau} \left( s + m
\, s_{\bar \tau}, \theta \right) + d^1_\tau \, \vec a + \delta_{n,
\tau} \, \vec e_2 +  \tilde w_\tau (s, \theta) \, \vec N^{\vec
a}_{\bar \tau} \left( s + m \, s_{\bar \tau}, \theta \right)
\]
where
\[
\tilde w_\tau (s, \theta ) : = \xi (s) \, w_\tau^1 \left( s + m \,
s_{\bar \tau}, \theta \right) + ( 1 - \xi (s)) \, w_\tau^1 \left(
m \, s_{\bar \tau} -s, \theta \right)
\]
We will denote by $A^1_{n,\tau}$ the image of $(-1,1 ) \times S^1$
by $Z^0_{n,\tau}$. We set
\[
Z^{\ell}_{n,\tau} : = {R_{2\pi\ell/k}} \, Z^0_{n,\tau}
\]
for $\ell =1, \ldots, k-1$ which describes the connection of the
image of $E^{1}_{n,\tau}$ by $R_{2\pi\ell/k}$ with the image of
$E^{-1}_{n,\tau}$ by $R_{2\pi (\ell+1)/k}$.

\subsection{Estimate of the mean curvature of the connected surface}

The compact surface which is obtained through these connections
will be denoted by $S_{n, \tau}$. It is an immersed compact
surface of genus $k$. By construction, the mean curvature of the
surface $S_{n, \tau}$ is equal to $1$ except in annular regions
$A_{n,\tau}^0$, $A_{n,\tau}^1$ and in their images by the elements
of $G_k$. The following estimates follow at once from the fact
that the functions $w^0_\tau$, $\bar w^0_\tau$ and $w^1_\tau$ are
exponentially decaying, as explained in \S 4.
\begin{lemma} We have
\[
\| H_{S_{n, \tau }} - 1 \|_{{\cal C}^{0,\alpha}(A^0_{n,\tau})}
\leq c \, e^{- n \, \gamma_{\tau ,2 } \, s_{\tau }}
\]
and, provided (\ref{eq:Zcondition}) is satisfied, we have
\[
\| H_{S_{n, \tau }} - 1 \|_{{\cal C}^{0,\alpha}(A^1_{n,\tau})}
\leq c \, e^{- m \, \gamma_{\bar \tau ,2} \, s_{\bar \tau}} .
\]
where the constant $c>0$ does not depend on $\tau \in I$ nor on $n
\in {\mathbb N}$. \label{lemma2}
\end{lemma}

\subsection{Partition of unity on $S_{n, \tau}$}

Subordinate to the above construction is a partition of unity we
now describe.

\begin{enumerate}
\item[(i)] The function $\chi_{n, \tau}$ is a smooth function
defined on $S_{n, \tau}$ which is equal to $1$ on
\[
\Sigma_{n,\tau} (n \, s_\tau-1, m \, s_{\bar \tau} -1) \subset
S_{n, \tau}
\]
and which is equal to $0$ on the complement of
\[
\begin{array}{rllll}
\Sigma_{n,\tau}(n \, s_\tau-1, m \, s_{\bar \tau} -1)  \cup
Y^0_{n,\tau} ((-1,1)\times S^1) & \cup&
Z^0_{n,\tau} ((-1,1)\times S^1) \\[3mm]
& \cup&  Z^{k-1}_{n,\tau} ((-1,1)\times S^1) \end{array}
\]
in $S_{n,\tau}$. To be more precise, on the part of $S_{n, \tau}$
parameterized by $Y^0_{n,\tau}$, the function $\chi_{n, \tau}$ is
equal to $1$ for $s \geq 1$ and equal to $0$ for $s \leq -1$ and
on the part of $S_{n, \tau}$ parameterized by $Z^0_{n,\tau}$, the
function $\chi_{n, \tau}$ is equal to $1$ for $s \leq -1$ and
equal to $0$ for $s \geq 1$. This function is assumed to be
invariant under the action of $S_1$.

\item[(ii)] The function $\bar \chi_{n, \tau}$ is a smooth
function defined on $S_{n, \tau}$ which is equal to $1$ on
\[
\bar \Sigma_\tau (n \, s_\tau -1) \subset S_{n, \tau}
\]
and which is equal to $0$ on the complement of
\[
\bar \Sigma_\tau (n \, s_\tau -1)  \cup_{\ell=0}^{k-1}
Y^\ell_{n,\tau} ((-1,1)\times S^1)
\]
in $S_{n, \tau}$. To be more precise, on the part of $S_{n, \tau}$
parameterized by $Y^0_{n,\tau}$, the function $\bar \chi_{n,
\tau}$ is equal to $1$ for $s \leq - 1$ and equal to $0$ for $s
\geq 1$. This function is assumed to be invariant under the action
of $G_k$.

\item[(iii)] We also ask that
\[
\bar \chi_{n, \tau} + \sum_{\ell =0}^{k-1}  \, \chi_{n, \tau}
\circ (R_{2\pi\ell/k})^{-1} = 1
\]
on $S_{n, \tau}$.
\end{enumerate}

There is another set of cutoff functions which will be needed.
They can be described as follows~:
\begin{enumerate}
\item[(i)] The function $\chi_{n, \tau}^e$ is a smooth function
defined on $S_{n, \tau}$ which is equal to $1$ on
\[
\begin{array}{rllll}
\Sigma_{n,\tau} (n \, s_\tau-1, m \, s_{\bar \tau}-1) \cup
Y^0_{n,\tau} ((- n \, s_\tau + 2 ,1)\times S^1) & \cup
Z^{0}_{n,\tau} ((-1, m\, s_{\bar \tau}-2)\times S^1) \\[3mm]
& \cup Z^{k-1}_{n,\tau} ((-1, m \, s_{\bar \tau}-2)\times S^1)
\end{array}
\]
and which is equal to $0$ on the complement of
\[
\begin{array}{rllll}
\Sigma_{n,\tau} (n \, s_\tau-1, m \, s_{\bar \tau}-1) \cup
Y^0_{n,\tau} ((- n \, s_\tau +1 ,1)\times S^1) & \cup
Z^0_{n,\tau} ((-1, m\, s_{\bar \tau}-1)\times S^1) \\[3mm]
& \cup  Z^{k-1}_{n,\tau} ((-1, m \, s_{\bar \tau}-1)\times S^1)
\end{array}
\]
To be more precise, on the part of $S_{n, \tau}$ parameterized by
$Y^0_{n,\tau}$, the function $\chi_{n, \tau}$ is equal to $1$ for
$s \geq - n \, s_\tau +2$ and equal to $0$ for $s \leq - n \,
s_\tau +1$ and on the part of $S_{n, \tau}$ parameterized by
$Z^0_{n,\tau}$, the function $\chi_{n, \tau}$ is equal to $1$ for
$s \leq  m \, s_{\bar \tau}-2$ and equal to $0$ for $s \geq
 m \, s_{\bar \tau} -1$. This function is assumed to be
invariant under the action of $S_1$.

\item[(ii)] The function $\bar \chi_{n, \tau}^e$ is a smooth
function which is equal to $1$ on
\[
\bar \Sigma_\tau (n\, s_\tau-1) \cup_{\ell=0}^{k-1}
Y^\ell_{n,\tau} ((0,  n  \, s_\tau -2)\times S^1)
\]
and which is equal to $0$ on the complement of
\[
\bar \Sigma_\tau (n\, s_\tau-1) \cup_{\ell=0}^{k-1}
Y^\ell_{n,\tau} ((0,  n  \, s_\tau -1)\times S^1)
\]
To be more precise, on the part of $S_{n, \tau}$ parameterized by
$Y^0_{n,\tau}$, the function $\chi_{n, \tau}$ is equal to $1$ for
$s \leq  n  \, s_\tau - 2$ and equal to $0$ for $s \geq  n \,
s_\tau -1$. This function is assumed to be invariant under the
action of $G_k$.
\end{enumerate}

\subsection{Extension of the elements of ${\mathcal K}_{
\Sigma_\tau }$ and ${\mathcal K}_{\bar \Sigma_\tau }$}

Building on the analysis of \S 4.3, we explain how the restriction
of the elements of ${\mathcal K}_{\Sigma_{n,\tau} }$ to
$\Sigma_{n,\tau} (n \, s_\tau -1, m\, s_\tau -1)\subset S_{n,
\tau}$ and the restriction of the elements of ${\mathcal K}_{\bar
\Sigma_\tau }$ to $\bar \Sigma_\tau (n\, s_\tau -1) \subset S_{n,
\tau}$ can be extended to functions which are defined on $S_{n ,
\tau }$. By "extension" we mean that we these restrictions are
first connected with restrictions of the elements of ${\mathcal
N}_{\Sigma_{n,\tau} }$ to $\Sigma_{n,\tau} (n \, s_\tau -1, m\,
s_\tau -1)$ and the restriction of the elements of ${\mathcal
N}_{\bar \Sigma_\tau }$  to $\bar \Sigma_\tau (n\, s_\tau -1)$ and
then extended to $S_{n,\tau}$ by using the action of $G_k$.  The
fact that these extensions are meaningful (see
Lemma~\ref{le:exto}) relies on (\ref{eq:Zcondition}). We keep the
same notations for the elements of ${\cal K}_{\Sigma_{n, \tau}}$
and ${\cal K}_{\Sigma_\tau}$.

\begin{itemize}
\item[(i)] The restriction of $\psi^{T,\vec a}_{\tau, E^1}$ and
$\psi^{T,\vec a^\perp}_{\tau, E^1}$ to $\Sigma_{n,\tau} (n\,
s_\tau, m \, s_{\bar \tau})$ can be easily extended to
$S_{n,\tau}$ using the fact that the  ends $E^1_\tau$ and
$R_{2\pi/k}\, E^1_\tau$ are symmetric with respect to $\Pi_k$. For
example, for $\vec b=\vec a,  \vec a^\perp$, we can first define a
function $\Psi^{T, \vec b}_{n, \tau}$ on the part of $S_{n, \tau}$
which is parameterized by $Z_{n,\tau}^1$,
\[
\Psi^{T, \vec b}_{ n, \tau}  = \chi_{n, \tau}  \, \psi^{T, \vec
b}_{\tau , E^1} + (1- \chi_{n,\tau} ) \,  \psi^{T, \vec b}_{\tau,
E^1} \circ ({R_{2\pi/k}})^{-1}
\]
and then use the action of $G_k$ to extend this function to the
other components of $S_{n, \tau}$.

\item[(ii)] The restriction of the element $\bar \psi^{T}_\tau$ of
${\mathcal K}_{\bar \Sigma_\tau}$ to $\bar \Sigma_\tau (n \,
s_\tau-1)$ can be extended to $S_{n,\tau}$ using the restriction
to $\Sigma_{n,\tau} (n, \, s_\tau -1, m, \, s_{\bar \tau}-1)$ of
$\Phi^{T, \vec e_2}_{\Sigma_{n,\tau}}$, the (unique) element of
${\mathcal N}_{\Sigma_{n,\tau}}$ which is asymptotic to $ \Phi^{T,
\vec e_2}_{E^0_{n,\tau}}$ on $E^0_{n,\tau}$ (i.e. the globally
defined Jacobi field which corresponds to translation of
$\Sigma_{n,\tau}$ along the $x_2$ axis). We define a function
$\bar \Psi^{T}_{n, \tau}$ first by writing
\[
\bar \Psi^{T}_{ n, \tau}  = \bar \chi_{n, \tau}  \, \bar
\psi^{T}_{\tau} + (1- \bar \chi_{n, \tau}) \, \Phi^{T, \vec
e_2}_{\Sigma_\tau}
\]
on the part of $S_{n, \tau}$ which is parameterized by $
Y_{n,\tau}^0$. Observe that $\Phi^{T, \vec e_2}_{\Sigma_\tau}$ is
asymptotic to a linear combination of $\psi^{T, \vec a}_{\tau,
E^1}$ and $\psi^{T, \vec a^\perp}_{\tau, E^1}$ on the other ends
of $\Sigma_\tau $ and we can use the extension described in (i) to
extend the function to $S_{n, \tau}$. For example,
\[
\bar \Psi^{T}_{n, \tau} = \chi_{n, \tau} \, \Phi^{T, \vec
e_2}_{\Sigma_\tau} + (1- \chi_{n,\tau} ) \, \Phi^{T, \vec
e_2}_{\Sigma_\tau} \circ ({R_{2\pi/k}})^{-1}
\]
on the part of $S_{n, \tau}$ which is parameterized by
$Z_{n,\tau}^0$, and then we use the action of $G_k$ to extend this
function to the other components of $S_{n, \tau}$.

\item[(iii)] One can choose the parameter $t$ in such a way that
the element $\psi^{D}_{\tau, E^0} + t \, \psi^T_{\tau, E^0}$ of
${\mathcal K}_{\Sigma_{n,\tau}}$ is asymptotic to $\Phi^{D}_{\bar
E^0_\tau}$, the generator of ${\mathcal N}_{\bar \Sigma_\tau}$.
Indeed, $ (\psi^{D}_{\tau, E^0} + t \, \psi^T_{\tau, E^0}) \circ
X^0_{ \tau}$ is asymptotic to $\Phi^{D}_{D_\tau} + t \, \Phi^{T,
\vec e_2}_{D_\tau}$ and  $\Phi^{D}_{\bar E^0_\tau} \circ \bar
X^0_\tau $ is asymptotic to $ \Phi^{D}_{D_\tau}$.  Granted the
definition of $Y^0_{n,\tau}$ (in terms of $\bar X^0_\tau$ and
$X^0_\tau$) together with (\ref{shift}), we choose
\[
t  =  n \, p_\tau \, s_\tau.
\]
These two functions are then connected, as in (i) or (ii), to
define the function $\Psi^D_{n, \tau}$. For example, we define
\[
\Psi^{D}_{ n, \tau}  = \bar \chi_{n, \tau}  \, \Phi^{D}_{\bar
E^0_\tau} + (1- \bar \chi_{n, \tau}) \, (\psi^{D}_{\tau, E^0} +
 n \, p_\tau \, s_\tau \, \psi^T_{\tau, E^0})
\]
on the part of $S_{n, \tau}$ which is parameterized by $
Y_{n,\tau}^0$ and then extend this function to all $S_{n, \tau}$
 using the action of $G_k$.
\end{itemize}

We define ${\cal L}_{S_{n,\tau}}$ to be the Jacobi operator about
the surface $S_{n, \tau}$. The following result again follows from
the fact that the functions $w^0_\tau$, $\bar w^0_\tau$ and
$w_\tau^1$ are exponentially decaying.
\begin{lemma}
There exists a constant $c>0$ which does not depend on $\tau \in
I$ nor on $n$ such that
\[
\|{\mathcal L}_{S_{n ,\tau}} \, \Psi^{T, \vec b}_{n,
\tau}\|_{{\mathcal C}^{0, \alpha}(A^1_{n,\tau})} \leq c \,  e^{-
\gamma_{\bar \tau , 2} \, m \, s_{\bar \tau}}
\]
for $\vec b = \vec a, \vec a^\perp$ and
\[
\|{\mathcal L}_{S_{n ,\tau}} \, \bar \Psi^{T}_{n, \tau}
\|_{{\mathcal C}^{0, \alpha}(A^0_{n,\tau})} \leq c  \, e^{-
\gamma_{\tau , 2} \, n \, s_\tau }  \qquad \mbox{and} \qquad \|
{\mathcal L}_{S_{n ,\tau}} \, \bar \Psi^{T}_{n, \tau}
\|_{{\mathcal C}^{0, \alpha}(A^1_{n,\tau})} \leq c \, e^{-
\gamma_{\bar \tau , 2} \, m \, s_{\bar \tau}}
\]
Finally, given $\mu \in (-\gamma_{\tau , 2} , 0)$, there exists a
constant $c_\mu>0$ which does not depend on $\tau \in I$ nor on
$n$ such that
\[
\| {\mathcal L}_{S_{n ,\tau}} \, \Psi^{D}_{n, \tau}\|_{{\mathcal
C}^{0, \alpha}(A^0_{n,\tau})}  \leq c_\mu \, e^{-\mu \, n
\,s_{\tau}}
\]
\label{le:exto}
\end{lemma}

\section{Perturbation of $S_{n, \tau}$}

\subsection{Mapping properties}

We define the weighted spaces on $S_{n,\tau}$. Roughly speaking,
to evaluate the norm in this space, we restrict a function to each
summand constituting $S_{n, \tau}$ and then evaluate each terms
using the norm defined in Definition~\ref{de:2}.
\begin{definition}
Given $r\in {\mathbb N}$, $\alpha\in(0,1)$ and $\mu \in {\R}$, we
define ${\cal C}^{r,\alpha}_\mu (S_{n, \tau })$ to be the space of
functions $ w \in  {\cal C}^{r,\alpha} (S_{n, \tau })$ which are
invariant under the action of $G_k$. This space is endowed with
the norm
\[
\| w \|_{{\cal C}^{r,\alpha}_\mu (S_{n, \tau } )} := \|
\chi_{n,\tau} \, w \|_{{\cal D}^{r,\alpha}_\mu (\Sigma_\tau )} +
\| \bar \chi_{n, \tau} \, w \|_{{\cal D}^{r,\alpha}_\mu (\bar
\Sigma_\tau)}
\]
\end{definition}

We also define the $4$ dimensional space
\[
{\mathcal K}_{ S_{n, \tau } } =  \mbox{Span} \{ \bar \Psi^{T}_{n,
\tau} , \Psi^{D}_{n, \tau} , \Psi^{T, \vec a}_{n, \tau} , \Psi^{T,
\vec a^\perp}_{n, \tau} \}
\]
In the following result we glue together the parametrices for
${\cal L}_{\Sigma_\tau }$ and ${\cal L}_{\bar \Sigma_\tau}$ to
obtain a parametrix for ${\cal L}_{S_{n,\tau}}$~:
\begin{proposition}
Assume that $\mu \in (-\inf (\gamma_{\tau ,2}, \gamma_{\bar \tau
,2}) , 0)$ is fixed. There exist $n_0>0$ and $c >0$ and, for all
$n \geq n_0$ and $\tau \in I$ for which (\ref{eq:Zcondition})
holds, one can find an operator
\[
G_{n, \tau} : {\cal C}_\mu^{0,\alpha}(S_{n, \tau })
\longrightarrow {\cal C}_\mu^{2,\alpha}(S_{n, \tau }) \oplus
{\calK}_{S_{n, \tau }},
\]
such that $w : = G_{n, \tau} (f)$ solves ${\calL}_{S_{n, \tau}} \,
w = f$ on $S_{n, \tau}$  and
\[
\| w \| _{{\cal C}_\mu^{2,\alpha}(S_{n, \tau }) \oplus
{\calK}_{S_{n, \tau }}} \leq c \, \| f \| _{{\cal
C}_\mu^{0,\alpha}(S_{n, \tau})},
\]
for some constant which does not depend on $\tau \in I$ nor on $n
\geq n_0$. \label{pr:ZAZA}
\end{proposition}
{\bf Proof~:} Given a function $g$ defined on $S_{n, \tau}$, it
will be convenient to identify the function $\chi_{n, \tau} \, g$
(resp. $\bar \chi_{n, \tau} \, g$) with a function which is
defined on $\Sigma_{n,\tau}$ (resp. $\bar \Sigma_\tau$). This
identification is done in the natural way on the common parts of
the surfaces and by identifying $(\chi_{n, \tau} \, g) \circ
Z^0_{n, \tau}$ with $(\chi_{n, \tau} g ) \circ X^1_{n,\tau}$,
$(\chi_{n, \tau} \, g) \circ Y^0_{n, \tau}$ with $(\chi_{n, \tau}
g ) \circ X^0_{n,\tau}$ and so on \ldots on the ends of the
surfaces.

\medskip

Conversely, given a function $g$ define in $\Sigma_{n,\tau}$
(resp. $\bar \Sigma_\tau$) we will identify the function $\chi_{n,
\tau}^e \, g$ (resp. $\bar \chi_{n, \tau}^e \, g$) with a function
which is defined on $\Sigma_\tau$ (resp. $\bar \Sigma_\tau$).

\medskip

Given $f\in {\cal C}_{\mu}^{0,\alpha}(S_{n, \tau}))$ we want to
solve the equation
\[
{\calL}_{S_{n, \tau}} \, w =  f
\]
on $S_{n, \tau}$. We solve
\[
{\cal L}_{\Sigma_{n,\tau}} \, w_1  = \chi_{n, \tau} \, f
\]
on $\Sigma_{n,\tau}$ and
\[
{\cal L}_{\bar \Sigma_\tau} \, w_2  = \bar \chi_{n, \tau} \, f
\]
on $\bar \Sigma_\tau $.

\medskip

The existence of $w_i$ follows at once from the analysis described
in \S 3 and we have the estimate
\begin{equation}
\|  w_1 \|_{{\cal D}_\mu^{2, \alpha}(\Sigma_{n,\tau} ) \oplus
{\cal K} (\Sigma_{n,\tau})} + \|  w_2 \|_{{\cal D}_\mu^{2,
\alpha}(\bar \Sigma_\tau ) \oplus {\cal K} (\bar \Sigma_\tau)}
\leq c \, \| f \| _{{\cal C}_\mu^{0,\alpha}(S_{n, \tau})}
\label{en3}
\end{equation}
where the constant $c >0$ does not depend on $n$ nor on $\tau\in
I$. Observe that the function $w_1$ can be decomposed as
\[
w_1:= v_1 + a_1 \, (\psi_{n, \tau}^{D} + t \, \psi^{T}_{n,\tau} )
+ b_1 \, \psi^{T, \vec a}_{n, \tau} + c_1 \,\psi^{T, \vec
a^\perp}_{n, \tau},
\]
and the function $w_2$ can be decomposed as
\[
w_2:= v_2 + a_2 \, \bar \psi^{T}_{n, \tau},
\]
This being understood, we define the function $w$ on $S_{n, \tau}$
by
\[
w = \chi^e_{\tau,n}\, v_1 + a_1 \, \Psi^{D}_{n, \tau} + b_1 \,
\Psi^{T, \vec a}_{n, \tau} + c_1 \,\Psi^{T,\vec a^\perp}_{n, \tau}
+ \bar \chi^e_{\tau ,n} \, v_2 + a_2 \, \bar \Psi_{n, \tau}^{T}.
\]
Observe that
\[
||w||_{{\cal C}^{2, \alpha}_\mu (S_{n, \tau}) \oplus {\cal
K}_{S_{n, \tau}}} \leq c \, ||f||_{{\cal C}^{0, \alpha}_\mu (S_{n,
\tau})}
\]
for some constant which does not depend on $n$ nor on $\tau\in I$.
We claim that
\[
||{\cal L}_{S_{n, \tau}} \, w - f ||_{{\cal C}^{0, \alpha}_\mu
(S_{n, \tau})} \leq c \, ( e^{ -\gamma_{\tau ,2} \, n \, s_\tau }
+ e^{2 \, \mu \, n \, s_\tau } + e^{ -\gamma_{\bar \tau ,2} \, m
\, s_{\bar \tau} } + e^{2 \, \mu \, m \, s_{\bar \tau} } ) \,  \|
f \|_{{\cal C}_\delta^{0,\alpha}(S_{n, \tau})}.
\]
Since our problem is invariant under the action of $G_k$, it is
enough to evaluate this quantity on $Y^0_{n, \tau} ((- n \,
s_{\tau}, n \, s_{\tau}) \times S^1)$ and on $Z^0_{n, \tau}((-m \,
s_{\bar \tau}, m\, s_{\bar\tau}) \times S^1)$. We focuss our
attention on the estimate of ${\cal L}_{S_{n, \tau}} \, w - f$ on
$Y^0_{n, \tau} ((- n \, s_{\tau},0) \times S^1)$, since the
estimates on the other parts can be obtained similarly.

\medskip

In $Y^0_{n, \tau} ((- n \, s_\tau + 2 , -1) \times S^1)$, all the
elements of ${\cal K}_{S_{n, \tau}}$ are pieces of Jacobi fields
in the sense that, for all $W \in {\cal K}_{S_{n, \tau}}$
\[
{\cal L}_{S_{n,\tau}} W =0
\]
in this set. Therefore,
\[
{\cal L}_{S_{n, \tau}} \, w - f = {\cal L}_{S_{n, \tau}} \,
(v_1+v_2) - f = {\cal L}_{S_{n, \tau}} \, v_1
\]
in this set since .

\medskip

We now use the fact that $Y^0_{n, \tau} ((- n \, s_\tau + 2 , -1)
\times S^1)$ can be considered as a normal graph over
$E^0_{n,\tau}$ for some function bounded  and whose derivatives
are bounded by a constant times $e^{-\gamma_{\tau ,2} (s+n \,
s_\tau)}$ in $(- n \, s_\tau , 0) \times S^1$. In particular, this
implies that
\[
{\cal L}_{S_{n, \tau}} - {\cal L}_{\Sigma_{n, \tau}}
\]
is a second order partial differential operator whose coefficients
as well as their derivatives are bounded by a constant times $e^{-
\gamma_{\tau , 2} (s + n \, s_\tau)}$ in $(- n \, s_\tau , 0)
\times S^1$. Since ${\cal L}_{\Sigma_{n, \tau}} \, v_1 =0$ in this
set, we conclude that
\[
\| e^{\mu s} ({\cal L}_{S_{n, \tau}} \, w - f) \|_{{\cal C}^{0,
\alpha} (Y^0_{n, \tau}((-n \, s_\tau+1, -1) \times S^1))} \leq c
\, (e^{2n\mu s_\tau} + e^{\gamma_{\tau ,2} \, n\, s_\tau}) \, \|
f\|_{{\cal C}_\delta^{0,\alpha}(S_{n, \tau})}
\]

In $Y^0_{n, \tau} ((- n \, s_\tau +1 , - n \, s_\tau +2) \times
S^1)$, we obtain, using similar arguments and taking into account
the influence of the cutoff function $\chi^e_{n, \tau}$
\[
\| e^{\mu s} ({\cal L}_{S_{n, \tau}} \, w - f) \|_{{\cal C}^{0,
\alpha} (Y^0_{n, \tau} ((-n \, s_\tau +1, -n \, s_\tau + 2) \times
S^1)) } \leq c\,
 e^{2n\mu s_\tau} \, \|
f\|_{{\cal C}_\delta^{0,\alpha}(S_{n, \tau})}
\]
and in $Y^0_{n, \tau} ((- 1  ,0 ) \times S^1)$, we obtain, using
similar arguments
\[
\| e^{\mu s} ({\cal L}_{S_{n, \tau}} \, w - f) \|_{{\cal C}^{0,
\alpha} (Y^0_{n, \tau} (-1 , 0) \times S^1))} \leq c\, - e^{-
\gamma_{\tau , 2} \, n\, s_\tau} \, \| f\|_{{\cal
C}_\delta^{0,\alpha}(S_{n, \tau})}
\]

\medskip

So far, we have produced a linear  operator
\[
\tilde G_{n, \tau} : {\cal C}_\mu^{0,\alpha}(S_{n, \tau })
\longrightarrow {\cal C}_\mu^{2,\alpha}(S_{n, \tau }) \oplus
{\calK}_{S_{n, \tau }},
\]
defined by $\tilde G_{n, \tau} (f) := w$, which is uniformly
bounded (with respect to $n \in {\mathbb N}$ and $\tau \in I$) and
which satisfies
\[
|||{\cal L}_{S_{n, \tau}} \circ \tilde G_{n, \tau} - I ||| \leq c
\,( e^{ -\gamma_{\tau , 2} \, n \, s_\tau } + e^{2 \, \mu \, n \,
s_\tau } + e^{ -\gamma_{\bar \tau , 2} \, m \, s_{\bar \tau} } +
e^{2 \, \mu \, m \, s_{\bar \tau} } ) .
\]
for some constant independent of $n \in {\N}$ and $\tau \in I$.
The result then follows from a simple perturbation argument,
provided $n$ is chosen large enough. \hfill $\Box$

\subsection{The nonlinear argument}

We define the functions
\[
\Lambda(\tau ):=  \frac{1}{T_{\bar \tau}} \left( \sin (\pi/k) \, (
d^0_\tau +  \bar d^0_\tau ) -  d^1_\tau \right)
\]
and
\[
\Gamma(\tau):= 2 \, \sin (\pi /k) \, \frac{T_{\tau}}{T_{\bar
\tau}}
\]
Recall that $\tau$ and $\bar \tau$ are related through
(\ref{eq:tta}). We now prove the main result of the paper~:
\begin{theorem}
There exists $n_0 > 0$ such that, for all $n \geq n_0$ and all
$\tau \in I$ satisfying
\begin{equation}
\Lambda(\tau) + n \, \Gamma(\tau) \in {\mathbb N}
\label{eq:xxx}
\end{equation}
the surface $\Sigma_{n,\tau}$ can be perturbed into
a constant mean curvature $1$ surface.
\end{theorem}
{\bf Proof~:} We consider surfaces which can be written as a
normal graph over $S_{n, \tau}$, for some function  $w \in {\cal
C}_\mu^{2,\alpha}(S_{n, \tau}) \oplus {\cal K}_{S_{n, \tau}}$. The
equation which guaranties that this surface has constant mean
curvature equal to $1$ can be written as
\begin{equation}
{\cal L}_{S_{n, \tau}} \,w + {\cal Q}_{n, \tau} (w) = 1 - H_{S_{n,
\tau}},
\label{end4}
\end{equation}
where ${\cal L}_{S_{n, \tau}}$ is the Jacobi operator about $S_{n,
\tau}$, $H_{S_{n, \tau}}$ is the mean curvature of $S_{n, \tau}$
and ${\cal Q}_{n, \tau}$ collects all the nonlinear terms. It
should be clear from the construction of $S_{n, \tau}$ that, given
$r \in {\mathbb N}$ there exists $c_r >0$ (independent of $\tau
\in I$ and of $n \in {\mathbb N}$) such that the following
pointwise bound holds
\begin{equation}
| {\cal Q}_{n, \tau} (w_2) - {\cal Q}_{n, \tau} (w_1) |_{{\cal
C}^r} \leq c_r \, ( |w_2 |_{{\cal C}^{r+2}}  + |w_1|_{{\cal
C}^{r+2}}  ) \, |w_2 -w_1|_{{\cal C}^{r+2}} \label{endo}
\end{equation}
provided $|w_1|_{{\cal C}^1} + |w_2|_{{\cal C}^1} \leq 1$, where
\[
|w |_{{\cal C}^r} = \sum_{j=0}^r |\nabla^j w|
\]
and partial derivatives are computed using the vector fields
$\del_s$ and $\del_\theta$ along the pieces of $S_{n,\tau}$
parameterized by $Y^\ell_{n, \tau}$ and $Z^\ell_{n, \tau}$ and
using a fixed set of vector fields (independent of $n$) away from
these pieces.

\medskip

We fix $\mu \in (-\inf(\gamma_{\tau , 2}, \gamma_{\bar \tau, 2})
,0)$. Using the result of Proposition~\ref{pr:ZAZA}, our problem
reduces to finding a fixed point for~:
\begin{equation}
F_{n, \tau} : w \longrightarrow  G_{n, \tau} \Bigr(1- H_{S_{n, r}}
- {\cal Q}_{n, \tau} (w)\Bigr). \label{map}
\end{equation}
which belongs to ${\cal C}_\mu^{2,\alpha}(S_{n, \tau}) \oplus
{\cal K}_{S_{n,\tau}}$. It follows from the result of
Lemma~\ref{lemma2} that
\[
|| 1 - H_{S_{n, \tau}} ||_{{\cal C}^{0, \alpha}_\mu (S_{n, \tau})}
\leq c \, ( e^{- (\gamma_{\tau, 2}  + \mu) \, n \, s_{\tau}} +
e^{- (\gamma_{\bar \tau ,2}  + \mu) \, m \, s_{\bar \tau}}).
\]
We set
\[
\rho_{n, \tau} : = ( e^{- (\gamma_{\tau ,2}  + \mu) \, n
\, s_{\tau}} + e^{- (\gamma_{\bar \tau , 2}  + \mu) \, m \,
s_{\bar \tau}}).
\]
Applying the result of Proposition~\ref{pr:ZAZA}, we conclude that
\begin{equation}
|| G_{n, \tau} (1 - H_{S_{n, \tau}}) ||_{{\cal C}^{2, \alpha}_\mu
(S_{n, \tau}) \oplus {\cal K}_{S_{n, \tau}} } \leq \bar c \,
\rho_{n, \tau}. \label{end1}
\end{equation}
for some constant $\bar c >0$ which does not depend on $\tau \in
I$ nor on $n\in {\mathbb N}$, for which (\ref{eq:xxx}) holds.

\medskip

Now, it follows from (\ref{endo}) that there exists a constant $c
>0$ which does not depend on $\tau\in I$ nor on $n \in {\mathbb
N}$ such that
\begin{equation}
|| {\cal Q}_{n, \tau} (w_2) - {\cal Q}_{n, \tau} (w_1) ||_{{\cal
C}^{0, \alpha}_\mu (S_{n, \tau}) } \leq c \, (n^2 + e^{ - \mu \, n
\, s_\tau} + e^{ - \mu \, m \, s_{\bar \tau}} ) \, \rho_{n, \tau}
\, ||w_2 -w_1||_{{\cal C}^{2, \alpha}_\mu (S_{n, \tau}) \otimes
{\cal K}_{S_{n,\tau}}} . \label{end}
\end{equation}
provided $\|w_2\|_{{\cal C}_\mu^{2,\alpha}(S_{n, \tau}) \oplus
{\cal K}_{S_{n,\tau}}} + \|w_1\|_{{\cal C}_\mu^{2,\alpha}(S_{n,
\tau}) \oplus {\cal K}_{S_{n,\tau}}} \leq 2 \, \bar c \, \rho_{n,
\tau}$. The $n^2$ which appears in this estimate arises from the
fact that the element $\Psi^D_{n, \tau}$ of ${\cal
K}_{S_{n,\tau}}$ is not bounded uniformly in $n$, but is bounded,
as well as its derivatives, by a constant (independent of $\tau$
and $n$) times $n$.

\medskip

We choose $\mu$  close enough to $0$ (but still negative !) so
that
\[
\lim_{n \rightarrow +\infty} ( e^{ - \mu \, n \, s_\tau} + e^{ -
\mu \, m \, s_{\bar \tau}} ) \, \rho_{n, \tau} =0
\]
uniformly for $\tau \in I$ (Recall that $n$ and $m$ are related by
(\ref{eq:xxx}), in particular there exists $c >0$, independent of
$\tau \in I$, such that $n \leq c\, m$ and $m\leq c\, n$). The
fact that, provided $n$ is chosen large enough, the mapping $F_{n,
\tau}$ has a fixed point in the ball of radius $2\, \bar c \,
\rho_{n, \tau}$ in ${\cal C}_\mu^{2,\alpha}(S_{n, \tau}) \oplus
{\cal K}_{S_{n,\tau}}$ follows directly from (\ref{end1}) and
(\ref{end}). \hfill $\Box$

\medskip

The surfaces we have obtained are immersed, compact surfaces with
genus $k$ (these surfaces are not embedded since the Type-$1$
elements which have been used for their construction are never
embedded). The surfaces obtained for different values of $\tau$
and $n$ satisfying (\ref{eq:xxx}) are geometrically different
(i.e. are not congruent modulo a rigid motion), provided $n_0$ is
chosen large enough. Hence, the set solutions of (\ref{eq:xxx})
give a local picture of the set of compact constant mean curvature
surfaces of genus $k$ with symmetry group $G_k$.

\medskip

Finally, observe that the result of Lemma~1 together with the fact
that $\tau$ and $\bar \tau$ are related by (\ref{eq:tta}) implies
that
\[
\partial_\tau \left( \frac{T_{\tau}}{T_{\bar \tau} } \right)  > 0
\]
Therefore (\ref{eq:xxx}) has nontrivial solutions $\tau \in I$,
for any $n$ large enough.

\end{document}